\documentstyle{amsppt}
\nologo \NoBlackBoxes \magnification1200 \pagewidth{6.5 true in}
\pageheight{9 true in} \topmatter
\title{On the number of linear forms in logarithms}
\endtitle
\author
Youness Lamzouri
\endauthor
\address{D{\'e}partment  de Math{\'e}matiques et Statistique,
Universit{\'e} de Montr{\'e}al, CP 6128 succ Centre-Ville,
Montr{\'e}al, QC  H3C 3J7, Canada}
\endaddress
\email{Lamzouri{\@}dms.umontreal.ca}
\endemail
\abstract Let $n$ be a positive integer.  In this paper we estimate
the size of the set of linear forms $b_1\log a_1 + b_2\log
a_2+...+b_n\log a_n$,
 where $|b_i|\leq B_i$ and $1\leq a_i\leq A_i$ are integers, as $A_i,B_i\to
 \infty$.
\endabstract
\endtopmatter

\document

\head 1. Introduction \endhead

 \noindent The theory of linear forms in logarithms, developed by A. Baker ([1] and [2]) in the 60's,
  is a powerful method in the
 transcendental number theory. It consists of finding lower bounds for $|b_1\log a_1 + b_2\log a_2+...+b_n\log a_n|$,
 where the $b_i$ are integers and the $a_i$ are algebraic numbers for which $\log a_i$
 are linearly independent over ${\Bbb Q}$. We consider the simpler case where the $a_i>0$ are integers, and
 we let $B_j=\max\{|b_j|,1\}$, and $B=\max_{1\leq j\leq n} B_j$.

 \noindent Lang and Waldschmidt [4]
  conjectured the following
\proclaim{Conjecture} Let $\epsilon>0$. There exists
$C(\epsilon)>0$ depending only on $\epsilon$, such that
$$ |b_1\log a_1 + b_2\log a_2+...+b_n\log
a_n|>\frac{C(\epsilon)^nB}{(B_1...B_na_1...a_n)^{1+\epsilon}}.
$$
\endproclaim
\noindent One part of the argument they used to motivate the
Conjecture, is that the number of distinct linear forms $b_1\log a_1
+ b_2\log a_2+...+b_n\log a_n$, where $|b_j|\leq B_j$ and $0<a_j\leq
A_j$, is $\asymp B_1...B_nA_1...A_n$, if the $A_i$ and the $B_i$ are
sufficiently large.

\noindent In this paper we estimate the number of these linear forms
as $A_i,B_i\to \infty$.

\noindent An equivalent formulation of the problem is to estimate
the size of the following set
$$ R=R(A_1,...,A_n,B_1,...,B_n):= \{r\in {\Bbb Q}: r=a_1^{b_1}a_2^{b_2}...a_n^{b_n},
 1\leq a_i\leq A_i, |b_i|\leq B_i \},$$ as $A_i,B_i \rightarrow
\infty.$

\noindent For the easier case $A_i=A$ and $B_i=B$ for all $i$, a
trivial upper bound on $|R|$ is

\noindent $2^nA^nB^n/n!+o(A^nB^n)$, since permuting the numbers
$a_i^{b_i}$ gives rise to the same number $r$.

\noindent We prove that this bound is attained asymptotically as
$A,B\rightarrow\infty$. Also we deal with the general case, which is
harder since not every permutation is allowed for all the ranges.
Indeed the size of $R$ depends on the ranges of the $A_i$ and the
$B_i$, as we shall see in Corollaries 1 and 2.
\smallskip

\noindent Let $E\subset\{(a_1,...,a_n,b_1,...,b_n),1\leq a_i\leq
A_i, |b_i|\leq B_i \}$. We say that $r\in {\Bbb Q}$ has a
representation in $E$, if $r=a_1^{b_1}a_2^{b_2}...a_n^{b_n}$, for
some $(a_1,...,a_n,b_1,...,b_n)\in E$.

\noindent  For  $r\in R$, if $\sigma\in S_n$ satisfies $1\leq
a_{\sigma(i)}\leq A_i$, and $| b_{\sigma(i)}|\leq B_i $ for all $i$,
we say that $\sigma$ permutes $r$, or $\sigma$ is a possible
permutation for the $a_i^{b_i}$. Finally we say that a permutation
$\sigma\in S_n$ is permissible if
$$ |\{ r\in R: \sigma \ \hbox{permutes} \
r\}| \gg A_1...A_nB_1...B_n.$$ The main result of this paper is the
following

\proclaim{Theorem } There exists a set
$E\subset\{(a_1,...,a_n,b_1,...,b_n),1\leq a_i\leq A_i, |b_i|\leq
B_i \}$ satisfying
$$ |E| \sim 2^n A_1A_2...A_nB_1B_2...B_n,$$ as
 $A_i,B_i\rightarrow \infty$, such that any
 rational number $r\in \{a_1^{b_1}a_2^{b_2}...a_n^{b_n} :(a_1,...,a_n,b_1,...,b_n)\in E \}$
 has a unique representation in $E$ up to permissible
 permutations.
 \endproclaim
 \noindent From this result we can deduce that $|R|$ is asymptotic to the cardinality of the set
 of $2n$-tuples
 $\{(a_1,...,a_n,b_1,...,b_n),1\leq a_i\leq A_i, |b_i|\leq
B_i \}$ modulo permissible permutations.

\noindent In the case $A_i=A$, $B_i=B$, every permutation is
permissible and we deduce the following Corollary

\proclaim{Corollary 1}As $A,B\rightarrow \infty$, we have
$$ |\{r\in {\Bbb Q}: r=a_1^{b_1}a_2^{b_2}...a_n^{b_n}, 1\leq a_i\leq A, |b_i|\leq B \}|=
\frac{2^nA^{n}B^n}{n!} + o(A^{n}B^n).$$
\endproclaim
\noindent Now suppose that $A_i=o(A_{i+1})$ for all $1\leq i\leq
n-1$, or $B_i=o(B_{i+1})$ for all $1\leq i\leq n-1$. For a
non-identity permutation $\sigma\in S_n$, there exists $j$ for which
$\sigma(j)\neq j$. Therefore if $\sigma$ permutes
$r=a_1^{b_1}a_2^{b_2}...a_n^{b_n}$, we must have $1\leq
a_j,a_{\sigma(j)} \leq \min(A_j,A_{\sigma(j)})$ and
$-\min(B_j,B_{\sigma(j)})\leq b_j,b_{\sigma(j)} \leq
\min(B_j,B_{\sigma(j)})$. And so we deduce that
 $$
 \align
 &|\{r\in {\Bbb Q}: r=a_1^{b_1}a_2^{b_2}...a_n^{b_n},
1\leq a_i\leq A_i, |b_i|\leq B_i : \sigma \ \hbox{permutes} \
r\}|\\
&\leq
2^nA_1...A_nB_1...B_n\left(\frac{\min(A_j,A_{\sigma(j)})\min(B_j,B_{\sigma(j)})}
{\max(A_j,A_{\sigma(j)})\max(B_j,B_{\sigma(j)})}\right)=o(A_1...A_nB_1...B_n),
\endalign
$$ by our assumption on the $A_i$ and $B_i$. Thus in this case no permutation $\sigma\neq 1$ is permissible.
Therefore we have

 \proclaim {Corollary 2} If
$A_i=o(A_{i+1})$ for all $1\leq i\leq n-1$, or $B_i=o(B_{i+1})$ for
all $1\leq i\leq n-1$, then
$$ |\{r\in {\Bbb Q}:  r=a_1^{b_1}a_2^{b_2}...a_n^{b_n},
 1\leq a_i \leq A_i, |b_i|\leq B_i
\}|\sim 2^nA_1...A_nB_1...B_n,$$ as $A_i,B_i\rightarrow \infty$.
\endproclaim

\noindent We can observe that Corollaries 1 and 2 correspond to
extreme cases: in Corollary 1 all permutations are permissible,
while none is permissible in Corollary 2. Indeed we can prove

\proclaim {Corollary 3} As $A_i,B_i\rightarrow \infty$, we have
$$ \frac{2^n}{n!}A_1...A_nB_1...B_n \lesssim |R| \lesssim
2^nA_1...A_nB_1...B_n.$$ Moreover the two bounds are optimal.
\endproclaim

\demo{ Proof } From the Theorem we have that
$$ |R| \sim \sum_{\Sb 1\leq a_1\leq A_1\\ |b_1|\leq B_1\endSb}...
\sum_{\Sb 1\leq a_n\leq A_n\\ |b_n|\leq B_n\endSb}
\frac{1}{|\{\sigma\in S_n: \sigma \hbox{ is possible for the }
a_i^{b_i} \}|}.$$ The result follows from the fact that $ 1\leq
|\{\sigma\in S_n: \sigma \hbox{ is possible for the } a_i^{b_i}
\}|\leq n!$.
\enddemo

\noindent For the simple case $n=2$, there is only one non-trivial
permutation $\sigma=(12)$. This permutation is possible only if
$1\leq a_1,a_2\leq \min(A_1,A_2)$ and $|b_1|,|b_2|\leq
\min(B_1,B_2)$. Then by the Theorem, and after a simple calculation
we deduce that
$$
\align &|\{r\in {\Bbb Q}: r=a_1^{b_1}a_2^{b_2}, 1\leq a_1\leq A_1,
1\leq a_2\leq A_2, |b_1|\leq B_1, |b_2|\leq B_2  \}|\\
&\sim 4A_1A_2B_1B_2- 2\min(A_1,A_2)^2\min(B_1,B_2)^2.
\endalign
$$
In general the size of $|R|$ is asymptotic to an homogeneous
polynomial of degree $2n$ in the variables
$A_1,...,A_n,B_1,...,B_n$.  Moreover it's also necessary to order
the $A_i$'s and $B_i$'s, so without loss of generality we assume
that $A_1\leq A_2\leq ...\leq A_n$ and $B_{\pi(1)}\leq
B_{\pi(2)}\leq ...\leq B_{\pi(n)}$, where $\pi\in S_n$ is a
permutation.

\noindent  We prove the following

\proclaim{ Proposition} Suppose that  $A_1\leq A_2\leq ...\leq A_n$
and $B_{\pi(1)}\leq B_{\pi(2)}\leq ...\leq B_{\pi(n)}$, where
$\pi\in S_n$ is a permutation. Also let $A_0=B_{\pi(0)}=1$.

\noindent Then $|R|$ is asymptotic to
$$ 2^n \sum_{\Sb i_1=1\\1\leq j_1\leq \pi^{-1}(1)\endSb}\sum_{\Sb 1\leq i_2\leq 2\\1\leq j_2\leq
\pi^{-1}(2)\endSb}...\sum_{\Sb 1\leq i_n\leq n\\1\leq j_n\leq
\pi^{-1}(n)\endSb}\frac{\prod_{k=1}^n(A_{i_k}-A_{i_k-1})(B_{\pi(j_k)}-B_{\pi(j_k-1)})}{|\{\sigma\in
S_n: i_{\sigma(l)}\leq l, j_{\sigma(l)}\leq \pi^{-1}(l), \forall
1\leq l\leq n\}|},$$ as $A_i,B_i\to\infty$.
\endproclaim

\noindent {\bf Acknowledgments.}

\noindent I sincerely thank my advisor, Professor Andrew Granville,
for suggesting the problem, for many valuable discussions, and for
his encouragement during the various stages of this work.
\newpage

 \head 2. Preliminary lemmas \endhead

\noindent Let $C$ be a positive real number.  We say that the
$n$-tuple $(a_1,...,a_n)$ satisfies condition ($1_C$), if there
exists a prime $p$, such that $p^k|a_1a_2...a_n$ where $k\geq 2$,
and $p^k\geq C$.
 \proclaim{ Lemma 1} We have
$$ |\{(a_1,...,a_n), 1\leq a_i\leq A_i: \hbox{which satisfy
($1_C$)}\}| \ll_n \frac{A_1...A_n(\log C)^{n}}{\sqrt{C}}.$$
\endproclaim

\demo{Proof} First we have
$$
\align &|\{(a_1,...,a_n), 1\leq a_i\leq A_i: \hbox{which satisfy
($1_C$)}\}| \\
 &\leq \sum_{p} |\{(a_1,...,a_n), 1\leq a_i\leq A_i :
\exists k\geq 2, p^k\geq C,  \hbox{and} \ p^k|a_1a_2...a_n \}|. \tag{1}\\
\endalign
$$
\noindent Case 1. $p\leq \sqrt{C}$

\noindent In this case pick $k$ to be the smallest integer such
that $p^k\geq C$, ie $k=[\log C/\log p]+1$. Then the number of
$(a_1,...,a_n)$ such that $p^k|a_1a_2...a_n$ is equal to
$$ \sum_{d_1d_2...d_n=p^k}\prod_{i=1}^n \sum\Sb 1\leq a_i\leq A_i \\
d_i|a_i\endSb 1 \leq d_n(p^k)\frac{A_1...A_n}{p^k}\leq
d_n(p^k)\frac{A_1...A_n}{C}.$$ Now $d_n(p^k)= \binom{n+k-1}{k}$,
and by Stirling's formula, for $k$ large enough we have
 $$
 \align
 \log d_n(p^k)&= (n+k-1+\frac{1}{2})\log(n+k-1)- (k+\frac{1}{2})\log k - (n-1+\frac{1}{2})\log(n-1)+O(1)\\
 &\leq (k+\frac{1}{2})\log\left(1+\frac{n-1}{k}\right) + (n-\frac{1}{2}) \log\left(\frac{n-1+k}{n-1}\right) \\
 & \leq n \log k . \\
 \endalign
 $$
Then summing over these primes gives
$$\sum_{p \leq \sqrt{C}} |\{(a_1,...,a_n), 1\leq a_i\leq A_i : p^k|a_1a_2...a_n \}| = O_n\left( \frac{A_1...A_n(\log
C)^n}{\sqrt{C}}\right).\tag{2}$$

\noindent Case 2. $p > \sqrt{C}$

\noindent In this case pick $k=2$. Then the number of
$(a_1,...,a_n)$ such that $p^2|a_1a_2...a_n$ is
$O(A_1...A_n/p^2)$, where the constant involved in the $O$ depends
only on $n$. Therefore summing over these primes gives
$$\sum_{p > \sqrt{C}}  |\{(a_1,...,a_n), 1\leq a_i\leq A_i : p^2|a_1a_2...a_n \}| = O_n\left(
\frac{A_1...A_n}{\sqrt{C}}\right). \tag{3}$$ Thus combining (1),
(2) and (3) gives the result.
\enddemo
\smallskip
\noindent We say that $(a_1,...,a_n)$ satisfies condition $(2_C)$ if
at least one of the $a_i$ is $C$-smooth: that is has all its prime
factors lying below $C$.

\proclaim{ Lemma 2} Write $C^{u_i}=A_i$ for all $1\leq i\leq n$.
Then uniformly for $\min_{1\leq i\leq n} A_i\geq C\geq 2$, we have

$$ |\{(a_1,...,a_n), 1\leq a_i\leq A_i : \hbox{which satisfy
($2_C$)}\}| \ll_n A_1A_2...A_n\left(\sum_{i=1}^n
e^{-u_i/2}\right).$$
\endproclaim
\demo{Proof} We have that
$$
|\{(a_1,...,a_n), 1\leq a_i\leq A_i : \hbox{which satisfy
($2_C$)}\}| \ll_n A_1A_2...A_n\sum_{i=1}^n
\frac{\Psi(A_i,C)}{A_i},$$ where $\Psi(x,y)$ is the number of
$y$-smooth positive integers below $x$. The result follows by the
following Theorem of de Bruijn [3]
$$ \Psi(A_i,C) \ll  A_i e^{-u_i/2},$$ uniformly for $A_i\geq C\geq
2$.
\enddemo

 \noindent We say that $(b_1,b_2,...,b_n)$ satisfy
condition $(3_C)$, if there exists an n-tuple of integers $|c_i|\leq
2\log C$ not all zero, such that $c_1b_1+c_2b_2+...+c_nb_n=0$.

\proclaim{Lemma 3} We have that
$$|\{ (b_1,...,b_n), |b_i|\leq B_i : \hbox{which satisfy
condition $(3_C)$}\}| \leq
B_1B_2...B_n\sum_{i=1}^n\left(\frac{(9\log C)^n}{B_i}\right).
$$
\endproclaim
\demo{Proof} We note that
$$
\align & |\{ (b_1,...,b_n), |b_i|\leq B_i : \hbox{which satisfy
condition $(3_C)$}\}| \\
& \leq \sum \Sb |c_i|\leq 2\log C \\
(c_1,...,c_n)\neq (0,...,0) \endSb |\{(b_1,...,b_n), |b_i|\leq B_i : c_1b_1+c_2b_2+...+c_nb_n=0\}|\\
&\leq \sum \Sb |c_i|\leq 2\log C \\
(c_1,...,c_n)\neq (0,...,0) \endSb
(2B_1+1)...(2B_n+1)\sum_{i=1}^n\left(\frac{1}{2B_i+1}\right)
\leq B_1B_2...B_n\sum_{i=1}^n\left(\frac{(9\log C)^n}{B_i}\right).\\
\endalign
 $$
 \enddemo
\head 3. Proof of the results \endhead
 \demo{Proof of the Theorem}
 We begin by choosing $C:= \min(B_1,...,B_n,\log A_1,...,\log A_n)$. We consider the following
 set
 $$
 \align
 E:=&\{ (a_1,...,a_n,b_1,...,b_n), 1\leq a_i\leq A_i, |b_i|\leq B_i :\\
  &\hbox{$(a_i)$
 don't satisfy any of $(1_C)$,$(2_C)$, $(b_i)$
 don't satisfy $(3_C)$}\}.\\
 \endalign
 $$
 Then by our choice of $C$, if we combine Lemmas 1, 2 and 3, we
 observe that

 \noindent $|E|=2^nA_1...A_nB_1...B_n(1+o(1))$.

 \noindent Therefore it remains to prove
 that any representation of a rational number $r$ as
 $a_1^{b_1}a_2^{b_2}...a_n^{b_n}$ where $(a_1,...,a_n,b_1,...,b_n)$ belongs to
  $E$, is unique up to possible permutations of the $a_i^{b_i}$, and
  finally we can consider only permissible permutations (since the
  number of $r\in R$ which can be permuted by a non-permissible
  permutation is negligible).

 \noindent We begin by considering the following equation
 $$a_1^{b_1}a_2^{b_2}...a_n^{b_n}=e_1^{f_1}e_2^{f_2}...e_n^{f_n},\tag{4}$$
 where $(a_1,...,a_n,b_1,...,b_n)$ and ($e_1,...,e_n,f_1,...,f_n)$ are in $E$. If for some $i$, $a_i$
 contains a prime factor $p$ such that $p^2\nmid a_1a_2...a_n$ and
 $p^2\nmid e_1e_2...e_n$, then $b_i\in
 \{f_1,f_2,...f_n\}$. Now suppose that there exists $1\leq j\leq n$ such that $b_j\notin
 \{f_1,f_2,...f_n\}$, then for all the primes $p$ that divides $a_j$, there exists $k\geq 2$
 for which
 $p^k|a_1a_2...a_n$ or $p^k|e_1e_2...e_n$, but the $(a_i)$ and the
 $(e_i)$ don't satisfy condition $(1_C)$ and so we must have
 $p^k\leq C$, which implies that $a_j$ is $C$-smooth; however this
 contradicts the fact that the $(a_i)$ do not satisfy condition
 $(2_C)$. Therefore we deduce that
$$ \{b_1,b_2,...,b_n\}=\{f_1,f_2,...,f_n\}.$$
Then up to permutations, we have that $b_i=f_i$, and so equation (4)
become
$$a_1^{b_1}a_2^{b_2}...a_n^{b_n}=e_1^{b_1}e_2^{b_2}...e_n^{b_n}.\tag{5}$$
Let $p$ be any prime dividing $a_1a_2...a_n$, and let $\alpha_i\geq
0$ and $\beta_i\geq 0$ be the corresponding powers of $p$ in $a_i$
and $e_i$ respectively, and let $c_i=\alpha_i-\beta_i$. Then
equation (5) implies that $$ c_1b_1+c_2b_2+...c_nb_n=0.$$ Now the
$(a_i)$ and the $(e_i)$ do not satisfy condition $(1_C)$, and so
$0\leq \alpha_i,\beta_i\leq \log C/\log2\leq 2\log C$, which implies
that $|c_i|\leq 2\log C$. And since the $(b_i)$ do not satisfy
condition $(3_C)$, we deduce that $c_i=0$, and then
$\alpha_i=\beta_i$ for all $1\leq i\leq n$. Since this is true for
every prime factor of $a_1a_2...a_n$, we must have $a_i=e_i$ for all
$1\leq i\leq n$, and our Theorem is proved.
\enddemo

\demo{ Proof of the Proposition } We want to count the number of
elements $r=(r_1,...,r_n)$, where $r_i=(a_i,b_i)\in
[1,A_i]\times[-B_i,B_i]\cap {\Bbb Z}\times {\Bbb Z}$, modulo
possible permutations of the $r_i$'s.

\noindent Since the number of $r$ for which some $b_i$ is $0$, is
$o(A_1...A_nB_1...B_n)$, we can suppose that all the $b_i$'s are
positive by symmetry.

\noindent Moreover let $R_i:=[1,A_i]\times[1,B_i]\cap {\Bbb Z}\times
{\Bbb Z}$, and define the following distinct discrete sets
$R_{ij}:=[A_{i-1},A_i]\times[B_{\pi(j-1)}, B_{\pi(j)}]\cap {\Bbb
Z}\times {\Bbb Z}$, for $1\leq i,j\leq n$.

\noindent For every $1\leq k\leq n$, we have
$$R_k=\bigsqcup_{\Sb 1\leq i_k\leq k\\
1\leq j_k\leq \pi^{-1}(k)\endSb}R_{i_kj_k}.\tag{6}$$ This implies
$$ R_1\times R_2\times...\times R_n= \bigsqcup_{\Sb i_1=1\\
1\leq j_1\leq \pi^{-1}(1)\endSb}\bigsqcup_{\Sb 1\leq i_2\leq 2\\
1\leq j_2\leq \pi^{-1}(2)\endSb}...\bigsqcup_{\Sb 1\leq i_n\leq n\\
1\leq j_n\leq \pi^{-1}(n)\endSb}R_{i_1j_1}\times
R_{i_2j_2}\times...\times R_{i_nj_n}.$$

\noindent Now consider the elements $r\in R_{i_1j_1}\times
R_{i_2j_2}...\times R_{i_nj_n}$, with $1\leq i_k\leq k$ and

\noindent $1\leq j_k\leq\pi^{-1}(k)$ being fixed. If  $\sigma\in
S_n$ permutes $r$, then $r_{\sigma(k)}\in R_k$ for all $1\leq k\leq
n$, but $r_{\sigma(k)} \in R_{i_{\sigma(k)}j_{\sigma(k)}}$ also,
which implies that $R_{i_{\sigma(k)}j_{\sigma(k)}}\bigcap R_k\neq
\emptyset$. From (6) this is equivalent to
$R_{i_{\sigma(k)}j_{\sigma(k)}}\subseteq R_k$, and thus to the fact
that $i_{\sigma(k)}\leq k$ and $j_{\sigma(k)}\leq \pi^{-1}(k)$ for
all $1\leq k\leq n$.

\noindent Therefore  for any $r\in R_{i_1j_1}\times
R_{i_2j_2}...\times R_{i_nj_n}$, the number of $\sigma\in S_n$ which
permutes $r$ is constant and equal to

$$|\{\sigma\in S_n: i_{\sigma(l)}\leq l,
j_{\sigma(l)}\leq \pi^{-1}(l), \forall 1\leq l\leq n\}|.$$

\noindent Thus the number of elements in $R_1\times R_2...\times
R_n$, modulo possible permutations is
$$ \sum_{\Sb i_1=1\\1\leq j_1\leq \pi^{-1}(1)\endSb}\sum_{\Sb 1\leq i_2\leq 2\\1\leq j_2\leq
\pi^{-1}(2)\endSb}...\sum_{\Sb 1\leq i_n\leq n\\1\leq j_n\leq
\pi^{-1}(n)\endSb}\frac{\prod_{k=1}^n(A_{i_k}-A_{i_k-1})(B_{\pi(j_k)}-B_{\pi(j_k-1)})}{|\{\sigma\in
S_n: i_{\sigma(l)}\leq l, j_{\sigma(l)}\leq \pi^{-1}(l), \forall
1\leq l\leq n\}|},$$ which implies the result.\enddemo

\enddocument